\documentclass[psamsfonts]{amsart}
\usepackage[utf8]{inputenc}
\usepackage{amsfonts}
\usepackage{hyperref}
\usepackage{amsmath}
\usepackage{xcolor}
\usepackage{amsthm}
\usepackage{pdflscape}
\usepackage{pgfplots}
\usepackage{mathrsfs}

\newtheorem{thm}{Theorem}[section]

\newtheorem{lem}[thm]{Lemma}

\newtheorem{ppty}[thm]{Property}

\theoremstyle{definition}
\newtheorem{defn}[thm]{Definition}

\theoremstyle{remark}

\makeatletter
\let\c@equation\c@thm
\makeatother
\numberwithin{equation}{section}

\title{On the quantum differentiation of smooth real-valued functions}

\author{Kolosov Petro}

\begin{document}

\begin{abstract}
Calculating the value of $C^{k\in\{1,\infty\}}$ class of smoothness real-valued function's derivative in point of $\mathbb{R}^+$ in radius of convergence of its Taylor polynomial (or series), applying an analog of Newton's binomial theorem and $q$-difference operator. $(P,q)$-power difference introduced in section 5. Additionally,  by means of Newton's interpolation formula, the discrete analog of Taylor series, interpolation using $q$-difference and $p,q$-power difference is shown.
\\\\
\noindent \textbf{Keywords.} derivative, differential calculus, differentiation, Taylor's theorem, Taylor's formula, Taylor's series, Taylor's polynomial, power function, Binomial theorem, smooth function, real calculus, Newton's interpolation formula, finite difference, q-derivative, Jackson derivative, q-calculus, quantum calculus, (p,q)-derivative, (p,q)-Taylor formula\\\\
\noindent \textbf{2010 Math. Subject Class.} 26A24, 05A30, 41A58\\
\noindent \textbf{e-mail:} \ kolosov\_94@mail.ua\\
\noindent \textbf{ORCID:} \ \href{http://orcid.org/0000-0002-6544-8880}{http://orcid.org/0000-0002-6544-8880} \\
\end{abstract}

\maketitle
\tableofcontents

\section{Introduction}
Let be Taylor's theorem (see \S7 "Taylor's formula", \cite{1})
\begin{thm}\label{ttt1}
\textbf{Taylor's theorem.} Let be $n\geq1$ an integer, let function $f(x)$ be $n+1$ times differentiable in neighborhood of $a\in\mathbb{R}$. Let $x$ be an any function's argument from such neighborhood, $p$ - some positive number. Then, there is exist some $c$ between points $a$ and $x$, such that
\begin{equation}\label{mvm}
  f(x)=f(a)+\frac{f'(a)}{1!}(x-a)+\frac{f''(a)}{2!}(x-a)^2+\cdots+\frac{f^{(n)}(a)}{n!}(x-a)^n+R_{n+1}(x)
\end{equation}
where $R_{n+1}(x)$ - general form of remainder term

\begin{equation}\label{xcxxxs}
R_{n+1}(x)=\left(\frac{x-a}{x-a}\right)^p \frac{(x-c)^{n+1}}{n!p}f^{(n+1)}(c)
\end{equation}
\begin{proof}
Denote $\varphi(x,a)$ polynomial related to $x$ of order $n$, from right part of (\hyperref[mvm]{1.2}), i.e
\begin{equation}\label{rer}
  \varphi(x,a)=f(a)+\frac{f'(a)}{1!}(x-a)+\frac{f''(a)}{2!}(x-a)^2+\cdots+\frac{f^{(n)}(a)}{n!}(x-a)^n
\end{equation}
Next, denote as $R_{n+1}(x)$ the difference
\begin{equation}\label{a1}
  R_{n+1}(x)=f(x)-\varphi(x,a)
\end{equation}
Theorem will be proven, if we will find that $R_{n+1}(x)$ is defined by (\hyperref[xcxxxs]{1.3}). Let fix some $x$ in neighborhood, mentioned in theorem \hyperref[ttt1]{1.1}. By definition, let be $x>a$. Denote by $t$ an variable, such that $t\in[a,x]$, and review auxiliary function $\psi(t)$ of the form
\begin{equation}\label{a2}
  \psi(t)=f(x)-\varphi(x,t)-(x-t)^pQ(x)
\end{equation}
where
\begin{equation}\label{a3}
  Q(x)=\frac{R_{n+1}(x)}{(x-a)^p}
\end{equation}
More detailed $\psi(t)$ could be written as
\begin{equation}\label{a4}
  \psi(t)=f(x)-f(t)-\frac{f'(t)}{1!}(x-t)-\frac{f''(t)}{2!}(x-t)^2-\cdots-\frac{f^{(n)}(t)}{n!}(x-t)^n
\end{equation}
$$-(x-t)^pQ(x)$$
Our aim - to express $Q(x)$, going from properties of introduced function $\psi(t)$. Let show that function $\psi(t)$ satisfies to all conditions of Rolle's theorem \cite{2} on $[a,x]$. From (\hyperref[a4]{1.8}) and conditions given to function $f(x)$, it's obvious, that function $\psi(t)$ continuous on $[a,x]$. Given $t=a$ in (\hyperref[a2]{1.6}) and keeping attention to equality (\hyperref[a3]{1.7}), we have
\begin{equation}\label{a5}
  \psi(a)=f(x)-\varphi(x,a)-R_{n+1}(x)
\end{equation}
Hence, by means of (\hyperref[a1]{1.5}) obtain $\psi(a)=0$. Equivalent $\psi(x)=0$ immediately follows from (\hyperref[a4]{1.8}). So, $\psi(t)$ on segment $[a,x]$ satisfies to all necessary conditions of Rolle's theorem \cite{2}. By Rolle's theorem, there is exist some $c\in[a,x]$, such that
\begin{equation}\label{a6}
  \psi'(c)=0
\end{equation}
Calculating derivative $\psi'(t)$, differentiating equality (\hyperref[a4]{1.8}), we have
\begin{equation}\label{a7}
 \psi'(t)=-f'(t)+\frac{f'(t)}{1!}-\frac{f''(t)}{2!}(x-t)+\frac{f''(t)}{2!}2(x-t)-\cdots
\end{equation}
$$+\frac{f^{(n)}(t)}{n!}n(x-t)^{n-1}-\frac{f^{(n+1)}(t)}{n!}(x-t)^n+p(x-t)^{p-1}Q(x)$$
It's seen that all terms in right part of (\hyperref[a7]{1.11}), except last two items, self-destructs. Hereby,
\begin{equation}\label{a8}
  \psi'(t)=-\frac{f^{(n+1)}(t)}{n!}(x-t)^n+p\cdot(x-t)^{p-1}Q(x)
\end{equation}
Given $t=c$ in (\hyperref[a8]{1.12}) and applying (\hyperref[a6]{1.10}), obtain
\begin{equation}\label{a9}
  Q(x)=\frac{(x-c)^{n-p+1}}{n!p}f^{(n+1)}(c)
\end{equation}
By means of (\hyperref[a9]{1.13}) and (\hyperref[a3]{1.7}), finally, we have
\begin{equation}\label{a10}
  R_{n+1}(x)=(x-a)^pQ(x)=\left(\frac{x-a}{x-a}\right)^p \frac{(x-c)^{n+1}}{n!p}f^{(n+1)}(c)
\end{equation}
Case $x<a$ is reviewed absolutely similarly. (see for reference \cite{1}, pp 246-247)
\\This proves the theorem.
\end{proof}
\end{thm}
Let function $f(x)\in C^{k}$ class of smoothness and satisfies to theorem \hyperref[ttt1]{1.1}, then its derivative by means of its Taylor's polynomial centered at $a\in\mathbb{R}$ in radius of convergence with $f(x)$ and linear nature of derivative, $(gf(x)+um(x))^{'}=gf^{'}(x)+um^{'}(x)$, is
\begin{equation}\label{juj1}
\frac{d}{dx}f(x)=\frac{f'(a)}{1!}\frac{d}{dx}(x-a)^{}+\frac{f''(a)}{2!}\frac{d}{dx}(x-a)^{2}+\cdots+\frac{f^{(k-1)}(a)}{(k-1)!}\frac{d}{dx}(x-a)^{k}
\end{equation}
$$+R_{k+1}'(x)$$
Otherwise, if $f\in C^{\infty}$ we have derivative of Taylor series \cite{5} of $f$ given the same conditions as (\hyperref[juj1]{1.15})
\begin{equation}\label{xzx1}
  \frac{d}{dx}f(x)=\frac{f'(a)}{1!}\frac{d}{dx}(x-a)^{}+\frac{f''(a)}{2!}\frac{d}{dx}(x-a)^{2}+\cdots+\frac{f^{(n)}(a)}{n!}\frac{d}{dx}(x-a)^{n}+
\end{equation}
$$+\cdots$$
Hence, derivative of function $f:1\leq C(f)\leq \infty$ \footnote{For example, let $f$ be a $k$-smooth function, then $C(f)=k$.} could be reached by differentiating of its Taylor's polynomial or series in radius of convergence, and consequently summation of power's derivatives being multiplied by coefficient, according theorem \hyperref[ttt1]{1.1}, over $k$ from $1$ to $t\leq\infty$, depending on class of smoothness. Hereby, the properties of power function’s differentiation holds, in particular, the derivative of power close related to Newton's binomial theorem \cite{4}.
\begin{lem}\label{lm2}
Derivative of power function equals to limit of Binomial expansion of $(x+\Delta x)^n$, iterated from $1$ to $n$, divided by $\Delta x:\Delta x\to 0$.
\begin{proof}
 \begin{equation}\label{mm1}
 \frac{d(x^n)}{dx}=\lim\limits_{\Delta x \to 0}\left\{\sum_{k=1}^n{n \choose k} x^{n-k} (\Delta x)^{k-1}\right\}={n \choose 1} x^{n-1}
 \end{equation}
\end{proof}
\end{lem}
According to lemma (\hyperref[lm2]{1.17}), Binomial expansion is used to reach derivative of power, otherwise, let be introduced expansion, based on  forward finite differences, discussed in \cite{3}
\begin{equation}\label{h2}
  x^{n}=x^{n-2}+j\sum\limits_{k\in\mathfrak{C}(x)}^{}k\cdot x^{n-2}-k^2\cdot x^{n-3}, \ \ \ \ x\in\mathbb{N}
\end{equation}
where $j=3! \ \mathrm{and} \ \mathfrak{C}(x):=\{0, \ 1, \ \ldots, \ x\}\subseteq \mathbb{N}$.
Particularize\footnote{Transferring $x^{n-2}$ under sigma operator, decreasing the power by 1 and taking summation over $k\in\mathfrak{U}(x)$} (\hyperref[h2]{1.19}), one has
\begin{equation}\label{h1}
  x^{n}=\sum\limits_{k\in\mathfrak{U}(x)}^{}j\cdot k\cdot x^{n-2}-j\cdot k^2\cdot x^{n-3}+x^{n-3}
\end{equation}
where $\mathfrak{U}(x):=\{0, \ 1, \ \ldots, \ x-1\}\subseteq \mathbb{N}$.
\begin{ppty}
Let $\mathfrak{S}(x)$ be a set $\mathfrak{S}(x):=\{1, \ 2, \ \ldots, \ x\}\subseteq \mathbb{N}$, let be (\hyperref[h1]{1.20}) written as $T(x, \ \mathfrak{U}(x))$, then we have equality
\begin{equation}\label{ppty11}
  T(x, \ \mathfrak{U}(x))\equiv T(x, \ \mathfrak{S}(x)), \ x\in \mathbb{N}
\end{equation}
Let (\hyperref[h2]{1.19}) be denoted as $U(x,\ \mathfrak{C}(x))$, then
\begin{equation}\label{ppty2}
  U(x,\ \mathfrak{C}(x))\equiv U(x,\ \mathfrak{S}(x)) \equiv U(x,\ \mathfrak{U}(x))
\end{equation}
\begin{proof}
Let be a plot of $jkx^{n-2}-jk^2x^{n-3}+x^{n-3}$ by $k$ over $\mathbb{R}_{\leq 10}^{+}$, given $x=10$
\begin{center}
\begin{tikzpicture}[scale=1]
\begin{axis}[
	xlabel = {$k\in\mathbb{R}_{\leq 10}$},
	ylabel = {$y=jkx^{n-2}-jk^2x^{n-3}+x^{n-3}$},
	minor tick num = 1
]
\addplot[smooth, blue] table [x = b, y = a] {
	  a        b
      1        0
      55       1
      97       2
      127      3
      145      4
      151      5
      145      6
      127      7
      97       8
      55       9
      1        10
};
\addplot[only marks, red] coordinates {
( 0, 1 )
( 1, 55 )
( 2, 97 )
( 3, 127 )
( 4, 145 )
( 5, 151 )
( 6, 145 )
( 7, 127 )
( 8, 97 )
( 9, 55 )
( 10, 1 )
};
\end{axis}
\end{tikzpicture}
\end{center}
\begin{center}
  Figure 1. Plot of $jkx^{n-2}-jk^2x^{n-3}+x^{n-3}$ by $k$ over $\mathbb{R}_{\leq 10}^{+}$, $x=10$
\end{center}
Obviously, being a parabolic function, it's symmetrical over $\tfrac{x}{2}$, hence equivalent $T(x, \ \mathfrak{U}(x))\equiv T(x, \ \mathfrak{S}(x)), \ x\in \mathbb{N}$ follows. Reviewing (\hyperref[h2]{1.19}) and denote $u(t)=tx^{n-2}-t^2x^{n-3}$, we can make conclusion, that $u(0)\equiv u(x)$, then equality of $U(x,\ \mathfrak{C}(x))\equiv U(x,\ \mathfrak{S}(x)) \equiv U(x,\ \mathfrak{U}(x))$ immediately follows.\\This completes the proof.
\end{proof}
\end{ppty}
By definition we will use set $\mathfrak{U}(x)\subseteq \mathbb{N}$ in our next expressions.\\
Since, for each $x=x_0\in\mathbb{N}$ we have equivalent
\begin{lem}\label{gggg2}
$\forall x=x_0\in\mathbb{N}$ holds
\begin{equation}\label{hh2}
\underbrace{\sum_{t=1}^{x}\sum_{k=1}^n{n \choose k} t^{n-k}}_{x^n}\equiv\sum\limits_{k=0}^{x-1}j\cdot k\cdot x^{n-2}-j\cdot k^2\cdot x^{n-3}+x^{n-3}
\end{equation}
\begin{proof}
Proof can be done by direct calculations.
\end{proof}
\end{lem}
By lemma \hyperref[gggg2]{1.24} we have right to substitute (\hyperref[h1]{1.20}) into limit (\hyperref[mm1]{1.18}), replacing Binomial expansion, and represent derivative of power by means of expression (\hyperref[h1]{1.20}). Note that,
\begin{equation}\label{kkk33}
  \Delta (x^n)=\sum_{k=1}^n{n \choose k} x^{n-k} \neq j\cdot k\cdot x^{n-2}-j\cdot k^2\cdot x^{n-3}+x^{n-3}
\end{equation}
As (\hyperref[h1]{1.20}) is analog of Binomial expansion of power and works only in space of natural numbers, different in sense, that Binomial expansion, for example, could be denoted as $M(x, \ \mathfrak{C}(n))$, where $n$ - exponent. While (\hyperref[h1]{1.20}) could be denoted $T(x, \ \mathfrak{U}(x)\equiv \mathfrak{S}(x))$, it shows that in case of Binomial expansion the set over which we take summation depends on exponent $n$ of initial function, when for (\hyperref[h1]{1.20}) it depends on point $x=x_0\in\mathbb{N}$.
To provide expressions' (\hyperref[h1]{1.20}) usefulness\footnote{By classical definition of derivative, we have to use upper summation bound $(x+\Delta x)\in\mathbb{R^+}$ on (\hyperref[h1]{1.20}), which turns false result as (\hyperref[h1]{1.20}) works in space of $\mathbb{N}$. } on taking power's derivative over $\mathbb{R}^+$, derivative in terms of quantum calculus should be applied, as next section dedicated to.
\section{Application of Q-Derivative}
Derivative of the function $f$ defined as limit of division of function's grow rate by argument's grow rate, when grow rate tends to zero, and graphically could be interpreted as follows
\begin{center}
\begin{tikzpicture}[scale=1] 
\draw[->] (0,0) -- (6,0) node[anchor=north] {$x$};
\draw[->] (0,0) -- (0,4) node[anchor=east] {$y$};
\draw[blue, thick] plot[smooth] coordinates {
(0.4,	0.36)
(0.8,	0.49)
(1.2,	0.64)
(1.6,	0.81)
(2,	1)
(2.4,	1.21)
(2.8,	1.44)
(3.2,	1.69)
(3.6, 1.96)
(4, 2.25)
(4.4, 2.56)
(4.8, 2.89)
};
\draw (-0.2,0) node[anchor=north west] {$0$};
\draw[densely dashed,thin] (1.57,0) -- (1.57, 0.8);
\draw (1.9,-0.1) node[anchor=north east] {$x_0$};
\draw[densely dashed,thin] (4,0) -- (4, 2.25);
\draw (4.82,-0.1) node[anchor=north east] {$x_0+\Delta x$};
\draw[densely dashed,thin] (0,0.8) -- (1.57, 0.8);
\draw[thin] (0.2,0.101910828) -- (5, 2.547770701);
\draw (10,2.9) node[anchor=north east] {tangent line at $x_0$ as $\Delta x \to 0^+$};
\draw (0,1.1) node[anchor=north east] {$f(x_0)$};
\draw (5.9,3.4) node[anchor=north east] {$f(x)$};
\draw[densely dashed,thin] (0,2.25) -- (4, 2.25);
\draw (0,2.6) node[anchor=north east] {$f(x_0+\Delta x)$};
\end{tikzpicture}
\end{center}
\begin{center}
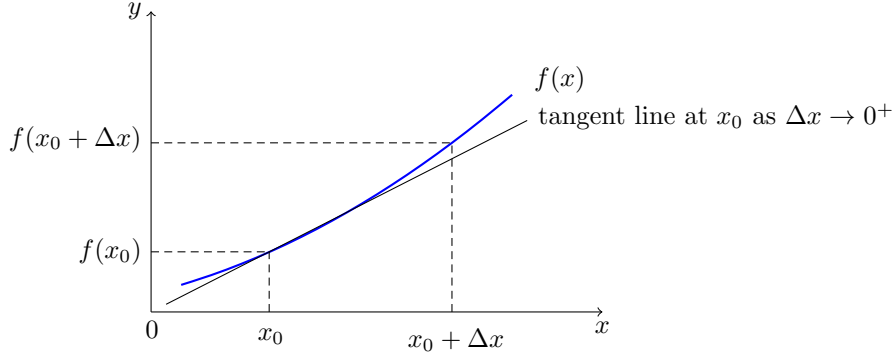
Figure 2. Geometrical sense of derivative
\end{center}
In 1908 Jackson \cite{10} reintroduced \cite{11}, \cite{12} the Euler-Jackson $q$-difference operator \cite{9}
\begin{equation}\label{x1}
 (D_qf)(x)=\frac{f(x)-f(qx)}{(1-q)x}, \ \  x\neq 0
\end{equation}
The limit as $q$ approaches $1^-$ is the derivative
\begin{equation}\label{x2}
  \frac{df}{dx}=\lim\limits_{q \to 1^-}(D_qf)(x)
\end{equation}
More generalized form of $q$-derivative
\begin{equation}\label{gdfgd22}
\frac{df(x)}{dx}=\lim\limits_{q \to 1^-}\underbrace{\frac{f(x)-f(xq)}{x-xq}}_{(D_qf^-)(x)}\equiv\lim\limits_{q \to 1^+} \underbrace{\frac{f(xq) -f(x)}{xq-x}}_{(D_qf^+)(x)}
\end{equation}
where $(D_qf^+)(x) \ \mathrm{and} \ (D_qf^-)(x)$ forward and backward $q$-differences, respectively.
The follow figure shows the geometrical sense of above equation as $q$ tends to $1^+$
\begin{center}
\begin{tikzpicture}[scale=1] 
\draw[->] (0,0) -- (6,0) node[anchor=north] {$x$};
\draw[->] (0,0) -- (0,4) node[anchor=east] {$y$};
\draw[blue, thick] plot[smooth] coordinates {
(0.4,	0.36)
(0.8,	0.49)
(1.2,	0.64)
(1.6,	0.81)
(2,	1)
(2.4,	1.21)
(2.8,	1.44)
(3.2,	1.69)
(3.6, 1.96)
(4, 2.25)
(4.4, 2.56)
(4.8, 2.89)
};
\draw (-0.2,0) node[anchor=north west] {$0$};
\draw[densely dashed,thin] (1.57,0) -- (1.57, 0.8);
\draw (1.9,-0.1) node[anchor=north east] {$x_0$};
\draw[densely dashed,thin] (4,0) -- (4, 2.25);
\draw (4.77,-0.1) node[anchor=north east] {$x_0\cdot q$};
\draw (5.9,3.4) node[anchor=north east] {$f(x)$};
\draw[densely dashed,thin] (0,0.8) -- (1.57, 0.8);
\draw[thin] (0.2,0.101910828) -- (5, 2.547770701);
\draw (10,2.9) node[anchor=north east] {tangent line at $x_0$ as $q \to 1^+$};
\draw (0,1.1) node[anchor=north east] {$f(x_0)$};
\draw[densely dashed,thin] (0,2.25) -- (4, 2.25);
\draw (0,2.6) node[anchor=north east] {$f(x_0\cdot q)$};
\end{tikzpicture}
\end{center}
\begin{center}\label{fig3}
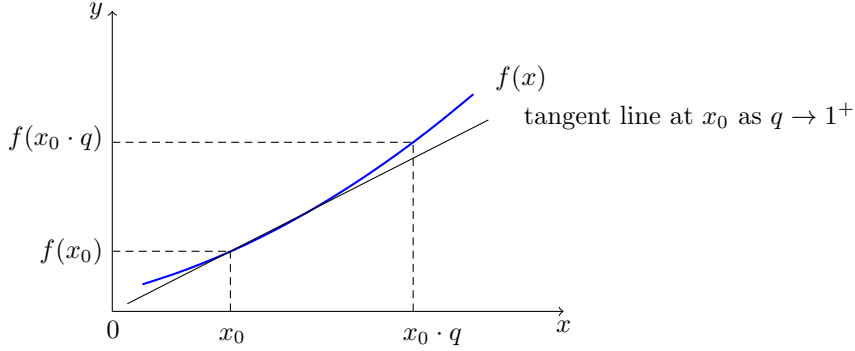
Figure 3. Geometrical sense of right part of (\hyperref[gdfgd22]{2.3})
\end{center}
Review the monomial $x^{n}$, where $n$-positive integer and applying right part of (\hyperref[gdfgd22]{2.3}), then in terms of $q$-calculus we have forward $q$-derivative over $\mathbb{R}$
\begin{equation}\label{sss22}
  \frac{d(x^n)}{dx}=\lim\limits_{q\to 1^+}(D_qx^{n^+})(x)=\lim\limits_{q\to 1^+}\frac{x^n(q^n-1)}{x(q-1)}
\end{equation}
$$=\lim\limits_{q\to 1^+}x^{n-1}\sum_{k=0}^{n-1}q^k, \ \ q\in\mathbb{R}$$
Otherwise, see reference \cite{9}, equation (109).\\
Generalized view of high-order power's derivative by means of (\hyperref[sss22]{2.4})
\begin{equation}\label{jkjk2}
  \frac{d^k(x^n)}{dx^k}=\lim\limits_{q\to 1^+}(D^k_qx^{n^+})(x)=\lim\limits_{q\to 1^+}x^{n-k}\prod_{j=0}^{k-1}\left(\sum_{m=0}^{n-j}q^m\right)
\end{equation}
Since, the main property of power is
\begin{ppty}\label{ppt1}
$$(x\cdot y)^n=x^n\cdot y^n$$
\end{ppty}
Let be definition
\begin{defn}\label{d2}
 By property (\hyperref[ppt1]{2.6}) and (\hyperref[h1]{1.20}), definition of $c=x\cdot t:t\in\mathbb{R}, \ x\in\mathbb{N}\Rightarrow c\in\mathbb{R}$ to power $n\in\mathbb{N}$
\begin{equation}\label{gfdgd2}
  c^n:=\xi(x,\ t)_n:=\sum\limits_{k=0}^{x-1}jkx^{n-2}\cdot t^n-jk^2x^{n-3}\cdot t^n+x^{n-3}\cdot t^n
\end{equation}
\end{defn}
Hereby, applying definition (\hyperref[d2]{2.7}) and (\hyperref[sss22]{2.4}), derivative of monomial $x^n:n\in\mathbb{N}$ by $x$ in point $x_0\in\mathbb{N}$ is\\
\begin{equation}\label{lll4}
  \frac{d(x^n)}{dx}\bigg|_{x=x_0}=\underbrace{\lim\limits_{q \to 1^+}\frac{\xi(x, \ q)_n-\xi(x, \ 1)_n}{x\cdot q-x}}_{\stackrel{\rm def}{=}\mathscr{D}_{q>1}[x^n]}\equiv\underbrace{\lim\limits_{q \to 1^-}\frac{\xi(x, \ 1)_n-\xi(x, \ q)_n}{x\cdot q-x}}_{\stackrel{\rm def}{=}\mathscr{D}_{q<1}[x^n]},
\end{equation}
Let us approach to extend the definition space of expression (\hyperref[lll4]{2.9}) from $x_0\in\mathbb{N}$ to $x_0\in\mathbb{R}^+$. Let be $x_0=\xi(t_0,\ p)_1\in\mathbb{R}^+\not\supseteq\mathbb{N}$ as $p\in\mathbb{R}^+\not\supseteq\mathbb{N}$ and $t_0\in\mathbb{N}$, then applying $(p,q)$-difference discussed in \cite{13}
\begin{equation}\label{xi}
  D_{p,q}f(x)=\frac{f(px)-f(qx)}{(p-q)x}, \  x\neq 0
\end{equation}
by means of definition (\hyperref[d2]{2.7}) and (\hyperref[xi]{2.10}), $(p,q)$-differentiating of monomial $x^n, \ n\in\mathbb{N}$ gives us
\begin{equation}\label{gfgf2}
\frac{d(x^n)}{dx}\bigg|_{x=t_0}=\lim\limits_{p \to q^+}D_{p,q}x^n=\underbrace{\lim\limits_{p \to q^+}\frac{\xi(x, \ p)_n-\xi(x, \ q)_n}{x\cdot p-x\cdot q}}_{\stackrel{\rm def}{=}\mathscr{D}_{p\rightarrow q}[x^n]}
\end{equation}
$$\equiv\underbrace{\lim\limits_{q \to p^-}\frac{\xi(x, \ p)_n-\xi(x, \ q)_n}{x\cdot p-x\cdot q}}_{\stackrel{\rm def}{=}\mathscr{D}_{p\leftarrow q}[x^n]},\ \ \ \ t_0\in\mathbb{N}, \ [p, \ q]\in\mathbb{R}^+\not\supseteq\mathbb{N}$$
Geometrical interpretation is shown below
\begin{center}
\begin{tikzpicture}[scale=1] 
\draw[->] (0,0) -- (6,0) node[anchor=north] {$x$};
\draw[->] (0,0) -- (0,4) node[anchor=east] {$y$};
\draw[blue, thick] plot[smooth] coordinates {
(0.4,	0.36)
(0.8,	0.49)
(1.2,	0.64)
(1.6,	0.81)
(2,	1)
(2.4,	1.21)
(2.8,	1.44)
(3.2,	1.69)
(3.6, 1.96)
(4, 2.25)
(4.4, 2.56)
(4.8, 2.89)
};
\draw (-0.2,0) node[anchor=north west] {$0$};
\draw[densely dashed,thin] (1.57,0) -- (1.57, 0.8);
\draw (1.9,-0.1) node[anchor=north east] {$x_0\cdot q$};
\draw[densely dashed,thin] (4,0) -- (4, 2.25);
\draw (4.82,-0.1) node[anchor=north east] {$x_0\cdot p$};
\draw[densely dashed,thin] (0,0.8) -- (1.57, 0.8);
\draw[thin] (0.2,0.101910828) -- (5, 2.547770701);
\draw (9.9,2.5) node[anchor=north east] {tangent line at $x_0$ as $p \to q^+$};
\draw (6.7,3.3) node[anchor=north east] {$f(x)=x^n$};
\draw (0,1.1) node[anchor=north east] {$\xi(x_0, \ q)_n$};
\draw[densely dashed,thin] (0,2.25) -- (4, 2.25);
\draw (0,2.6) node[anchor=north east] {$\xi(x_0, \ p)_n$};
\end{tikzpicture}
\end{center}
\begin{center}
 Figure 4. Geometrical interpretation of (\hyperref[gfgf2]{2.11})
\end{center}
\section{Application on functions of finite class of smoothness}
In this section we will get derivative of function $f\in C^n$ in point $x_0\in\mathbb{R}^+$ by means of its Taylor's polynomial and (\hyperref[lll4]{2.9}), where $n$ - some positive integer. Let $f(x)$ be an $n$-smooth function, then derivative of its Taylor's polynomial at radius of convergence with $f$ in $x_0:(x_0-a)\in\mathbb{N}$ is
\begin{equation}\label{dd}
  \frac{df(x)}{dx}\bigg|_{x=x_0}=\sum_{k=1}^{n}\left[\frac{f^{(k)}(a)}{k!}\mathscr{D}_{q>1}[(x-a)^k]\right]+\mathscr{D}_{q>1}[R_{n+1}(x)]
\end{equation}
$$\equiv\sum_{k=1}^{n}\left[\frac{f^{(k)}(a)}{k!}\mathscr{D}_{q<1}[(x-a)^k]\right]+\mathscr{D}_{q<1}[R_{n+1}(x)]$$
Otherwise, let $(x_0-a)$ satisfies to conditions of (\hyperref[gfgf2]{2.11}), i.e $(x_0-a)\in\mathbb{R}^+$, then applying operator $\mathscr{D}$, defined in (\hyperref[lll4]{2.9}) we can reach derivative of $f:f\in C^n$ in point $x_0:(x_0-a)\in\mathbb{R}^+$, by differentiation of its Taylor's polynomial in radius of convergence with $f$, that is
\begin{equation}\label{ddd}
  \frac{df(x)}{dx}\bigg|_{x=t_0}=\sum_{k=1}^{n}\left[\frac{f^{(k)}(a)}{k!}\mathscr{D}_{p\rightarrow q}[(x-a)^k]\right]+\mathscr{D}_{p\rightarrow q}[R_{n+1}(x)]
\end{equation}
$$\equiv\sum_{k=1}^{n}\left[\frac{f^{(k)}(a)}{k!}\mathscr{D}_{p\leftarrow q}[(x-a)^k]\right]+\mathscr{D}_{p\leftarrow q}[R_{n+1}(x)]$$
\section{Application on analytic functions}
If $f\in C^{\infty}$ (i.e analytic), then approximation by means of Taylor series holds in neighborhood of its center at $a\in\mathbb{R}$. Suppose that $f$ is real-valued and satisfies to conditions of Taylor's theorem \hyperref[ttt1]{1.1}, then derivative of $f$ at $x_0:x<x_0<a$ is
\begin{equation}\label{gbgb1}
  \frac{df(x)}{dx}=\frac{d}{dx}\sum_{k=0}^{\infty}\frac{f^{(k)}(a)}{k!}(x-a)^k=\left[\sum_{k=0}^{\infty}\frac{f^{(k)}(a)}{k!}\frac{d}{dx}(x-a)^k\right]_{x=x_0}
\end{equation}
Let $x_0$ satisfies to conditions of (\hyperref[dd]{3.1}), then, applying definition (\hyperref[d2]{2.7}), we have derivative of $f$ in point $x_0\in\mathbb{R}^+$
\begin{equation}\label{vvv2}
  \frac{df(x)}{dx}=\left[\sum_{k=1}^{\infty}\frac{f^{(k)}(a)}{k!}\mathscr{D}_{q>1}[(x-a)^k]\equiv\sum_{k=1}^{\infty}\frac{f^{(k)}(a)}{k!}\mathscr{D}_{q<1}[(x-a)^k]\right]_{x=x_0}
\end{equation}
Otherwise, if $x_0$ satisfies to conditions of (\hyperref[ddd]{3.2}) and $x_0$ in radius of convergence with $f$, then derivative of $f\in C^{\infty}$, by means of its Taylor's series and (\hyperref[d2]{2.7}), is
\begin{equation}\label{nmnm2}
 \frac{df(x)}{dx}=\left[\sum_{k=1}^{\infty}\frac{f^{(k)}(a)}{k!}\mathscr{D}_{p\rightarrow q}[(x-a)^k]\equiv\sum_{k=1}^{\infty}\frac{f^{(k)}(a)}{k!}\mathscr{D}_{p\leftarrow q}[(x-a)^k]\right]_{x=t_0}
\end{equation}
\section{Introduction of $(P, \ q)$-power difference}
\begin{lem}\label{lemmm2}
Let be $m\in\mathbb{R}/\mathbb{I}$ and $m$ could be represented as $m=at$, then exists some $c\in\mathbb{R}/\mathbb{I}$, such that
\begin{equation}
  m=a^c
\end{equation}
\end{lem}
Reviewing (\hyperref[gdfgd22]{2.3}), we can see, that argument's differential $\Delta x$ is given by $x\cdot q-x$, according to lemma \hyperref[lemmm2]{5.1} $\exists c\in\mathbb{R}/\mathbb{I}, \ x\cdot t-x=x^c-x$, then, from (\hyperref[sss22]{2.4}) immediately follows $q$-power difference, (see \cite{14}, page 2, equation 3)
\begin{equation}\label{xpxp}
  \mathcal{D}_{q>1}f(x):=\frac{f(x^q)-f(x^1)}{x^q-x^1}, \ \ x\neq 0
\end{equation}
As $q$ tends to $1^+$ we have reached derivative
\begin{equation}\label{xzxs}
  \frac{df(x)}{dx}=\lim\limits_{q\to 1^+}\mathcal{D}_{q>1}f(x)=\underbrace{\lim\limits_{q\to 1^+}\frac{f(x^q)-f(x^1)}{x^q-x^1}}_{\stackrel{\rm def}{=}\mathbf{D}_{q>1}[f(x)]}
\end{equation}
$$\equiv\underbrace{\lim\limits_{q\to 1^-}\frac{f(x^1)-f(x^q)}{x^1-x^q}}_{\stackrel{\rm def}{=}\mathbf{D}_{q<1}[f(x)]}=:\lim\limits_{q\to 1^-}\mathcal{D}_{q<1}f(x)$$
where $\lim\limits_{q\to 1^-}\mathcal{D}_{q<1}f(x)$ denotes the derivative through backward $q$-power difference
By lemma \hyperref[lemmm2]{5.1} from (\hyperref[xi]{2.10}) immediately follows $(p,q)$-power difference
\begin{equation}\label{pp1}
\mathcal{D}_{p\rightarrow q}f(x):=\frac{f(x^p)-f(x^q)}{x^p-x^q}, \ \ x\neq 0
\end{equation}
Hence, for $v=x^p, \ p\in\mathbb{R}$
\begin{equation}\label{bvbs2}
  \frac{df(x)}{dx}(v)=\lim\limits_{p\to q^+}\mathcal{D}_{p\rightarrow q}f(x)=\underbrace{\lim\limits_{p\to q^+}\frac{f(x^p)-f(x^q)}{x^p-x^q}}_{\stackrel{\rm def}{=}\mathbf{D}_{p\rightarrow q}[f(x)]}
\end{equation}
$$\equiv\underbrace{\lim\limits_{q\to p^-}\frac{f(x^p)-f(x^q)}{x^q-x^p}}_{\stackrel{\rm def}{=}\mathbf{D}_{p\leftarrow q}[f(x)]}=:\lim\limits_{q\to p^-}\mathcal{D}_{p\leftarrow q}f(x)$$
where $\mathbf{D}_{p\rightarrow q}[f(x)], \ \mathbf{D}_{p\leftarrow q}[f(x)]$ denote derivative through forward and backward $(p,q)$-power differences. Let us to show geometrical interpretation of (\hyperref[xzxs]{5.4}) and (\hyperref[bvbs2]{5.6})
\begin{center}
\begin{tikzpicture}[scale=1] 
\draw[->] (0,0) -- (6,0) node[anchor=north] {$x$};
\draw[->] (0,0) -- (0,4) node[anchor=east] {$y$};
\draw[blue, thick] plot[smooth] coordinates {
(0.4,	0.36)
(0.8,	0.49)
(1.2,	0.64)
(1.6,	0.81)
(2,	1)
(2.4,	1.21)
(2.8,	1.44)
(3.2,	1.69)
(3.6, 1.96)
(4, 2.25)
(4.4, 2.56)
(4.8, 2.89)
};
\draw (-0.2,0) node[anchor=north west] {$0$};
\draw[densely dashed,thin] (1.57,0) -- (1.57, 0.8);
\draw (1.9,-0.1) node[anchor=north east] {$x_0$};
\draw[densely dashed,thin] (4,0) -- (4, 2.25);
\draw (4.4,-0.1) node[anchor=north east] {$x_0^q$};
\draw (5.9,3.4) node[anchor=north east] {$f(x)$};
\draw[densely dashed,thin] (0,0.8) -- (1.57, 0.8);
\draw[thin] (0.2,0.101910828) -- (5, 2.547770701);
\draw (9.8,2.9) node[anchor=north east] {tangent line at $x_0$ as $q \to 1^+$};
\draw (0,1.1) node[anchor=north east] {$f(x_0)$};
\draw[densely dashed,thin] (0,2.25) -- (4, 2.25);
\draw (0,2.6) node[anchor=north east] {$f(x_0^q)$};
\end{tikzpicture}
\\Figure 5. Geometrical sense of (\hyperref[xzxs]{5.4})
\end{center}
\begin{center}
\begin{tikzpicture}[scale=1] 
\draw[->] (0,0) -- (6,0) node[anchor=north] {$x$};
\draw[->] (0,0) -- (0,4) node[anchor=east] {$y$};
\draw[blue, thick] plot[smooth] coordinates {
(0.4,	0.36)
(0.8,	0.49)
(1.2,	0.64)
(1.6,	0.81)
(2,	1)
(2.4,	1.21)
(2.8,	1.44)
(3.2,	1.69)
(3.6, 1.96)
(4, 2.25)
(4.4, 2.56)
(4.8, 2.89)
};
\draw (-0.2,0) node[anchor=north west] {$0$};
\draw[densely dashed,thin] (1.57,0) -- (1.57, 0.8);
\draw (1.9,-0.1) node[anchor=north east] {$x_0^q$};
\draw[densely dashed,thin] (4,0) -- (4, 2.25);
\draw (4.4,-0.1) node[anchor=north east] {$x_0^p$};
\draw (5.9,3.4) node[anchor=north east] {$f(x)$};
\draw[densely dashed,thin] (0,0.8) -- (1.57, 0.8);
\draw[thin] (0.2,0.101910828) -- (5, 2.547770701);
\draw (9.8,2.9) node[anchor=north east] {tangent line at $x_0$ as $p \to q^+$};
\draw (0,1.1) node[anchor=north east] {$f(x_0^q)$};
\draw[densely dashed,thin] (0,2.25) -- (4, 2.25);
\draw (0,2.6) node[anchor=north east] {$f(x_0^p)$};
\end{tikzpicture}
\\Figure 6. Geometrical sense of (\hyperref[bvbs2]{5.6})
\end{center}
Applying (\hyperref[xzxs]{5.4}) with monomial $x^m:m\in\mathbb{N}$, we get
\begin{equation}\label{bbb2}
\frac{d(x^m)}{dx}=\mathbf{D}_{q>1}[x^m]=\lim\limits_{q\to 1^+}\left[\sum_{k=1}^{m}(x^q)^{m-k}\cdot x^{k-1}\right]=mx^{m-1}
\end{equation}
$$\equiv\lim\limits_{q\to 1^-}\mathbf{D}_{q<1}[x^m]=\lim\limits_{q\to 1^-}\left[\sum_{k=1}^{m}x^{k-1}\cdot(x^q)^{m-k}\right]=mx^{m-1}$$
Note that $\mathbf{D}_{q<1}[x^m], \ \mathbf{D}_{q>1}[x^m]$ defined by (\hyperref[xzxs]{5.4}). The high order $N\leq m$ derivative, derived from (\hyperref[bbb2]{5.7})
\begin{equation}\label{bbb3}
\frac{d^N(x^m)}{dx^N}=\mathbf{D}_{q>1}^N[x^m]=\lim\limits_{q\to 1^+}\prod_{j=0}^{N-1}\left(\sum_{k=1}^{m-j}(x^q)^{m-k}\cdot x^{k-j-1}\right)
\end{equation}
$$\equiv\mathbf{D}_{q<1}^N[x^m]=\lim\limits_{q\to 1^-}\prod_{j=0}^{N-1}\left(\sum_{k=1}^{m-j}x^{k-j-1}\cdot(x^q)^{m-k}\right)$$
Let be analytic function $f$ and let $f$ satisfies to Taylor's theorem \hyperref[ttt1]{1.1} on segment of $(a,x), \ a\in\mathbb{R}$, then, applying (\hyperref[xzxs]{5.4}), in radius of convergence of its Taylor's series, we obtain derivative
\begin{equation}\label{lml2}
  \frac{df(x)}{dx}=\sum_{k=1}^{\infty}\frac{f^{(k)}(a)}{k!}\mathbf{D}_{q<1}[(x-a)^k]\equiv\sum_{k=1}^{\infty}\frac{f^{(k)}(a)}{k!}\mathbf{D}_{q>1}[(x-a)^k]
\end{equation}
Using $\mathbf{D}_{p\rightarrow q}[f(x)], \ \mathbf{D}_{p\leftarrow q}[f(x)]$ defined by  (\hyperref[bvbs2]{5.8}), for each $v=x^{p}$, we receive
\begin{equation}\label{pqp2}
\frac{df(x)}{dx}=\left[\sum_{k=1}^{\infty}\frac{f^{(k)}(a)}{k!}\mathbf{D}_{p\rightarrow q}[(x-a)^k]\equiv\sum_{k=1}^{\infty}\frac{f^{(k)}(a)}{k!}\mathbf{D}_{p\leftarrow q}[(x-a)^k]\right]_{x=v}
\end{equation}
Or, by means of definition (\hyperref[d2]{2.7}) and (\hyperref[bbb2]{5.9}), when $(x_0-a)\in\mathbb{N}$ derivative could be taken as follows
\begin{equation}\label{lml3}
\frac{df(x)}{dx}=\sum_{k=1}^{\infty}\left\{\frac{f^{(k)}(a)}{k!}\cdot\lim\limits_{n\to 1^+}\sum_{k=1}^{m}\xi(x-a, \ 1)_{nm-nk}x'\cdot \xi(x-a, \ 1)_{k-1}x'\right\}\Bigg|_{x=x_0}
\end{equation}
Given $x_0$, such that $(x_0-a)\in\mathbb{R}^+$, then conditions of (\hyperref[ddd]{3.2}) is reached, and, applying definition (\hyperref[d2]{2.7}), derivative $f'$ follows
\begin{equation}\label{lml5}
\frac{df(x)}{dx}=\sum_{k=1}^{\infty}\left\{\frac{f^{(k)}(a)}{k!}\cdot\lim\limits_{n\to 1^+}\sum_{k=1}^{m}\xi(x-a, \ 1)_{nm-nk}x'\cdot \xi(x-a, \ 1)_{k-1}x'\right\}\Bigg|_{x=t_0}
\end{equation}
\\Otherwise, let be $f:f\in C^{n}$, where $n$ - positive integer, then under similar conditions as (\hyperref[lml2]{5.11}) and (\hyperref[lml3]{5.13}), derivative could be reached by differentiating of $n$-order Taylor's polynomial of $f$ in terms of $q$-power difference (\hyperref[xpxp]{5.3}) under limit notation over $n$
\begin{equation}\label{lml4}
  \frac{df(x)}{dx}=\sum_{k=1}^{n}\left\{\frac{f^{(k)}(a)}{k!}\cdot\lim\limits_{n\to 1^+}\sum_{k=1}^{m}(x-a)^{nm-nk}x'\cdot (x-a)^{k-1}x'\right\}+R_{n+1}'(x)
\end{equation}
Similarly, as (\hyperref[lml3]{5.13}), derivative of $f\in C^n$ in point $x=x_0$, such that $(x_0-a)\in\mathbb{N}$
\begin{equation}\label{lml6}
\frac{df(x)}{dx}=\sum_{k=1}^{n}\left\{\frac{f^{(k)}(a)}{k!}\cdot\lim\limits_{n\to 1^+}\sum_{k=1}^{m}\xi(x-a, \ 1)_{nm-nk}x'\cdot \xi(x-a, \ 1)_{k-1}x'\right\}\Bigg|_{x=x_0}+R_{n+1}'(x)
\end{equation}
Otherwise, going from (\hyperref[lml5]{5.14}), $\forall (x_0-a)\in\mathbb{R}^+$
\begin{equation}\label{lml6}
\frac{df(x)}{dx}=\sum_{k=1}^{n}\left\{\frac{f^{(k)}(a)}{k!}\cdot\lim\limits_{n\to 1^+}\sum_{k=1}^{m}\xi(x-a, \ 1)_{nm-nk}x'\cdot \xi(x-a, \ 1)_{k-1}x'\right\}\Bigg|_{x=t_0}+R_{n+1}'(x)
\end{equation}
\section{Newton's interpolation formula}
Being a discrete analog of Taylor's series, the Newton's interpolation formula \cite{6}, first published in his Principia Mathematica in 1687, hereby, by author's opinion, supposed to be discussed
\begin{equation}\label{n1}
  f(x)=\sum_{k=0}^{\infty}\binom{x-a}{k}\Delta^k f(a)
\end{equation}
Given $q=const$ in (\hyperref[gdfgd22]{2.3}) divided $q$-difference $f[xq; \ x]$ is reached. Let be $\Delta f=f[xq; \ x](xq-x)$, then, by means of generalized high order forward finite difference $\Delta^k f, \ k\geq2$, (\cite{7}, \cite{8}), revised according to (\hyperref[gdfgd22]{2.3}), Newton's formula (\hyperref[n1]{6.1}) takes the form
\begin{equation}\label{n2}
  f(x)=\sum_{k=0}^{\infty}\left[\binom{x-a}{k}\sum_{m=0}^{k}(-1)^m\binom{m}{k}f(x\cdot t^m)\right]
\end{equation}
Review (\hyperref[xzxs]{5.4}) and given $q=const$ divided $q$-power difference follows, by similar way as (\hyperref[n2]{6.2}) reached, (\hyperref[n1]{6.1}) could be written as
\begin{equation}\label{n3}
  f(x)=\sum_{k=0}^{\infty}\left[\binom{x-a}{k}\sum_{m=0}^{k}(-1)^m\binom{m}{k}f(x^{n^{k-m}})\right]
\end{equation}
\section{Conclusion}
In this paper was discussed a way of obtaining real-valued smooth function's derivative in radius of convergence of it's Taylor's series or polynomial by means of analog of Newton's binomial theorem (\hyperref[h1]{1.20}) in terms of $q$-difference (\hyperref[dd]{3.1}) and $(p,q)$-power difference operators (\hyperref[pqp2]{5.12}). In the last section reviewed a discrete analog of Taylor's series - Newton's interpolation formula (\hyperref[n1]{6.1}), and applying operators of $q$-difference, $(p,q)$-power difference interpolation of initial function is shown (\hyperref[n2]{6.2}), (\hyperref[n3]{6.3}).

\end{document}